\documentclass{article}

\usepackage{amssymb}
\usepackage{amsfonts}
\usepackage{mathrsfs}
\usepackage{amsmath}
\usepackage{latexsym}
\usepackage{color}
\usepackage[utf8]{inputenc}

\newcommand{\be}{\begin{equation}}
\newcommand{\ee}{\end{equation}}
\newcommand{\bea}{\begin{eqnarray}}
\newcommand{\eea}{\end{eqnarray}}
\newcommand{\bean}{\begin{eqnarray*}}
\newcommand{\eean}{\end{eqnarray*}}
\newcommand{\la}{\label}

\newtheorem{theorem}{Theorem}
\newtheorem{lemma}{Lemma}
\newtheorem{coro}{Corollary}

 \newcommand{\nl}{\newline}

\newcommand{\R}{{\mathbb{R}}}
\newcommand{\C}{{\mathbb{C}}}

 \newcommand{\cE}{{\cal E}}

\newcommand{\po}{\partial\Omega}

 \newcommand{\as}{\mbox{ as }}

 \newcommand{\inprod}[2]{{\langle{#1},{#2}\rangle}}

\newcommand{\ia}{({\rm i})}
\newcommand{\ib}{({\rm ii})}

\newcommand{\ab}{\alpha\beta}
\newcommand{\lf}{\lambda\phi}
\newcommand{\subab}{\substack{|\alpha|=m \\ |\beta|=m}}
\newcommand{\subabl}{\substack{|\alpha|\leq m \\ |\beta|\leq m}}
\newcommand{\beq}{\begin{equation}}
\newcommand{\eeq}{\end{equation}}
\newcommand{\dom}{{\rm Dom}}
\newcommand{\comb}[2]{\big(\substack{ {#1} \\ {#2}}\big)}

\title{Gaussian estimates with best constants for higher-order Schr\"{o}dinger
operators with Kato potentials}
\date{}

\author{
G. Barbatis\footnote{Department of Mathematics,
 National and Kapodistrian University of Athens,  15784 Athens, Greece}}

\begin{document}

\date{\today}

\maketitle

\begin{abstract}
We establish Gaussian estimates on the heat kernel of a higher-order uniformly elliptic Schr\"{o}dinger operator 
with variable highest order coefficients and with a Kato class potential. The estimates involve the sharp constant in the Gaussian exponent.
\end{abstract}

\vspace{11pt}


\vspace{2pt}

 \parindent=25pt
 \parskip=0pt

\section{Introduction}

Let $\Omega$ be a domain in $\R^n$ and let $H_0$ be a uniformly elliptic operator of order $2m$ with $L^{\infty}$ coefficients acting on $L^2(\Omega)$,
\[
(H_0u)(x)=(-1)^m\sum_{\subabl}D^{\alpha}\{a_{\ab}(x)D^{\beta}u\},
\]
subject to Dirichlet boundary conditions on $\po$. If $2m\geq n$ then the semigroup generated by $H_0$ has a continuous integral
kernel $K_0(t,x,y)$ (also referred to as the heat kernel) which satisfies a Gaussian estimate of the form
\begin{equation}
\label{1}
| K_0(t,x,y)| < c_1t^{-\frac{n}{2m}}\exp\Big\{ -c_2\frac{|x-y|^{\frac{2m}{2m-1}}}{t^{\frac{1}{2m-1}}} + c_3t \Big\},
\end{equation}
see \cite{davies3,bg}.

In the article \cite{dh} Davies and Hinz studied the operator $H_0+V$ for singular potentials $V$ and obtained conditions under which 
the $L^2$ semigroup $e^{-(H_0+V)t}$ extends to a strongly continuous holomorphic semigroup in $L^p$, $1\leq p<\infty$. Amongst the potentials they considered are potentials $V$ that are
Kato class with respect to $H_0$, that is they satisfy
\[
  \| V(H_0+ \lambda)^{-1} \|_{L^1 \to L^1} \longrightarrow 0  \;\; , \qquad \as \lambda \to +\infty \, ,
\]
In the recent article \cite{ddy} the authors consider the question of Gaussian heat kernel estimates for $H_0+V$ for
Kato class potentials $V$. Under the assumption that $H_0$ has constant coefficients they prove that estimate (\ref{1}) is also valid for the heat kernel $K(t,x,y)$ of $H_0+V$. In the very recent article \cite{uwzd} the authors consider the operator $(-\Delta)^m +V$
for Kato potentials $V$ and apply the methods of \cite{ddy} together with Davies' exponential perturbation technique
as adpted in \cite{bd} in order to obtain estimates such as (\ref{1}) for $K(t,x,y)$ with the sharp constant $c_2$ in the Gaussian exponent.

The purpose of the present note is to show that if $2m>n$ then more can be achieved by an adaptation of the methods of \cite{b2001}.
Using purely $L^2$ methods we obtain a sharp Gaussian estimate for the heat kernel of $H_0+V$ for operators $H_0$ with variable coefficients. Moreover, unlike the three above mentioned articles, the Kato condition is imposed only on the negative part
$V_-$ of $V$, the positive part $V_+$ being merely in $L^1_{loc}(\Omega)$. 

The sharpness of these estimates depends of course on using the right distance function which is not the Euclidean but, rather, a Finsler distance induced by the operator. The sharp constant $\sigma_m$, also obtained in \cite{uwzd}, was first identified by Evgrafov and Postnikov \cite{ep} who obtained short time asymptotics of $K_0(t,x,y)$ for operators with constant coefficients in $\R^n$ and so-called {\em strongly convex} principal symbol (see definition below).

We prove two theorems which differ on the regularity assumptions imposed on the coefficients.
Theorem \ref{thm1} applies to operators with strongly convex symbol and  coefficients that
are bounded in the H\"{o}lder class $C^{m-2,1}(\Omega)$. Theorem \ref{thm2} is a more general result were the coefficients are merely in $L^{\infty}(\Omega)$ and the symbol need not be strongly convex;
the price paid is that instead of the sharp constant $\sigma_m$ we now have a constant $\sigma=\sigma_m -D$, with $D$ a certain measure of regularity for $H$.

\section{Setting and statement of results}\label{section:hke}

Let $\Omega$ be a domain in $\R^n$ and let $H_0$ be a uniformly elliptic operator of order $2m$ acting on $L^2(\Omega)$,
\[
(H_0u)(x)=(-1)^m\sum_{\subabl}D^{\alpha}\{a_{\ab}(x)D^{\beta}u\}
\]
subject to Dirichlet boundary condintions on $\po$. 
The coefficients $a_{\ab}(x)$, $|\alpha|,|\beta|\leq m$, are assumed to be real-valued functions in $L^{\infty}(\Omega)$ and the
matrix $\{a_{\ab}(x)\}$ is assumed to be symmetric for a.e. $x\in\Omega$ (the conditions on lower-order coefficients can easily be weakened). Under these assumptions the quadratic form
\[
Q_0(u)=\int_{\Omega}\sum_{\subabl}a_{\ab}(x)D^{\alpha}u D^{\beta}\bar{u}\, dx
\]
is then defined on $\dom(Q_0):= H^m_0(\Omega)$; we assume that G{\aa}rding's inequality
\begin{equation}
Q_0(u)\geq c_1\|u\|_{H^m(\Omega)}^2 -c_2 \|u\|_{L^2(\Omega)}^2  \; , \qquad u\in H^m_0(\Omega),
\label{barb:garding}
\end{equation}
is satisfied for some $c_1,c_2>0$. The quadratic form $Q_0(\cdot)$ is then closed and $H_0$ is defined on $L^2(\Omega)$ as the self-adjoint operator  associated to $Q_0(\cdot)$.
We note \cite[Theorem 7.12]{agmon} that inequality (\ref{barb:garding}) implies that the principal symbol of $H_0$ satisfies
\[
\sum_{\subab}a_{\ab}(x)\xi^{\alpha+\beta} \geq c_1|\xi|^{2m}, \; \qquad x\in\Omega , \; \xi\in\R^n .
\]

It is proved in \cite{davies1} that if $2m > n$ then the semigroup $e^{-H_0t}$ has a continuous integral kernel $K_0(t,x,y)$ which satisfies (\ref{1}). This result was later extended in the case $2m=n$ \cite{ElRo,AuTc}.
Estimate (\ref{1}) implies that the semigroup $e^{-H_0z}$, ${\rm Re}\, z>0$, extends to a strongly continuous bounded holomorphic semigroup $T_p(z)$ on $L^p(\Omega)$ for all $1\leq p <\infty$, and moreover the corresponding generators have spectrum which is independent of $p$ \cite{davies1}. In the case $2m<n$ the estimate (\ref{1}) is not valid as is seen by the counterexamples constructed in \cite{davies2}. We refer to the recent review article \cite{bg} for more information.

\subsection{Finsler distance and strong convexity}

To state our results we need to to define the distance function in terms of which our Gaussian estimates will be expressed
and also to introduce the notion of strong convexity.

The principal symbol
\[
A(x,\xi)=\sum_{\subab} a_{\ab}(x)\xi^{\alpha+\beta} \;\; , \quad\quad x\in\Omega \, , \; \xi\in\R^n,
\]
of $H$ induces canonically a Finsler distance $d(\cdot,\cdot)$ on $\Omega$ given by
\begin{equation}
d(y_1,y_2) =\sup\{ \phi(y_2) -\phi(y_1) \; : \; \mbox{ $\phi$ Lipschitz in $\Omega$ and
$A(x,\nabla\phi(x))\leq 1$ a.e. $x\in\Omega$}\}.
\label{finslermetric}
\end{equation}
If additional regularity is imposed on the coefficients then $d(\cdot,\cdot)$ is the distance induced by the Finsler metric with length element
$ds=p(x,dx)$ where
\[
p(x,\eta)= \sup_{\substack{\xi\in\R^n \\ \xi\neq 0}}\frac{\inprod{\xi}{\eta}}{A(x,\xi)^{1/2m}}\; , \quad x\in\Omega, \; \eta\in\R^n.
\]
This metric is Riemannian if $m=1$ or, more generally, if $A(x,\xi)$ is the $m$th power of a second order polynomial in $\xi$;
we refer to \cite{agmon1,barbatis1999} for a very short introduction to Finsler geometry and to \cite{bcz} for further reading.

Let the functions $a_{\gamma}$, $|\gamma|=2m$, be defined by
\[
A(x,\xi)=\sum_{|\gamma|=2m}\comb{2m}{\gamma} a_{\gamma}(x)\xi^{\gamma} \;\; , \qquad x\in\Omega \, , \; \xi\in\R^n ,
\]
where $\comb{2m}{\gamma} =(2m)!/ ( \gamma_1 !  \ldots \gamma_n !)$. The following notion of
{\em strong convexity} was first introduced by Evgrafov and Postnikov \cite{ep}.

\

\noindent
{\bf Definition.} The principal symbol $A(x,\xi)$ is strongly convex if for a.e. $x\in\Omega$ the quadratic form
\[
\Gamma(x ; p)=\sum_{\subab}a_{\alpha +\beta}(x)p_{\alpha}\overline{p_{\beta}} 
\]
is positive semi-definite on $\oplus_{|\alpha|=m}\C$.

Evgrafov and Postnikov \cite{ep} proved that if in addition to the assumptions above $H_0$ has constant coefficients on $\R^n$ and if the symbol $A(\xi)$ is strongly convex then
\[
K_0(t,x,y)  = \exp\Big\{ -\sigma_m \frac{d(x,y)^{2m/(2m-1)}}{t^{1/(2m-1)}} \Big\}(1+o(1)) , \qquad \quad \mbox{ as }t\to 0+ \, ,
\]
modulo subexponential terms, where
\[
\sigma_m =(2m-1)(2m)^{-\frac{2m}{2m-1}} \sin\big(\frac{\pi}{4m-2}\big).
\]
This was generalized to operators with smooth coefficients by Tintarev \cite{tintarev}.

The Gaussian estimates of Theorems \ref{thm1} and \ref{thm2} are expressed not in terms of $d(\cdot,\cdot)$ but rather in terms of an approximating family of distances: for any $M>0$ we define the distance
\[
d_M(y_1,y_2)=\sup\{ \phi(y_2)-\phi(y_1)  : \; \phi\in\cE_{A,M}\},
\]
where
\[
\cE_{A,M} =\big\{  \phi\in C^m(\Omega)  \; : \; A(x,\nabla\phi (x)) \leq 1  \; , \;
| D^{\alpha}\phi(x)| \leq M \; , \; x\in\Omega, \; 2\leq |\alpha| \leq m\big\}.
\]

\subsection{Kato potentials}\label{sec:kp}

Let $H_0$ be an operator of order $2m>n$ as above. We consider a real potential $V=V_+-V_-$ ($V_{\pm}\geq 0$)
and we make the following

\

\noindent
{\bf Hypothesis (H)} {\em The potentials $V_{\pm}$ belong in $L^1_{loc}(\Omega)$.
Moreover $V_-$  has zero form bound with respect to $H_0$, that is for any $\epsilon >0$ there exists
$c_{\epsilon}$ such that}
\beq
\int_{\Omega}V_- |u|^2dx \leq \epsilon Q_0(u) + c_{\epsilon}\|u\|_{L^2}^2 \; , \qquad u\in H^m_0(\Omega).
\label{ma}
\eeq
Under Hypothesis (H) the operator $H=H_0+V$ is defined in a standard way by means of the quadratic form
\begin{equation}
Q(u) =Q_0(u) +  \int_{\Omega}V|u|^2 dx  ,
\label{q}
\end{equation}
defined initially in $C^{\infty}_c(\Omega)$ and then extended by closure.
We note that (\ref{ma}) implies
\beq
Q_0(u) \leq \frac{1}{1-\epsilon}\Big(Q(u)+c_{\epsilon}\|u\|_{L^2}^2\Big) , \qquad u\in \dom(Q),
\label{q_0q}
\eeq
for any $\epsilon\in (0,1)$.

\

\noindent
{\bf Example 1.} If $V_-$ is Kato class with respect to $H_0$, that is if
\begin{equation}
\lim_{\lambda\to +\infty} \| V_-(H_0+\lambda)^{-1}\|_{L^1 \to L^1} =0 \, ,
\label{l1}
\end{equation}
then Hypothesis (H) is satisfied. This well known fact is seen by considering the weighted $L^p$ spaces
$L^p_{V_-}:= L^p(\Omega, V_-dx)$ (the fact that $V_-$ may be zero on a set of positive measure can easily be dealt with). We then have
\begin{eqnarray*}
\|  ( H_0 + \lambda)^{-1}V_-\|_{L^1_{V_-} \to L^1_{V_-}} &=&
 \sup_{w\in L^1_{V_-} } \frac{ \int_{\Omega}V_-(x)  \Big| \{ ( H_0 + \lambda)^{-1}(V_-w)\}(x) \Big| dx}{ \int_{\Omega}V_-(x)|w(x)|dx} \\
 &=&\sup_{w\in L^1_{V_-} } \frac{ \int_{\Omega}V_-(x)  \Big| \{ ( H_0 + \lambda)^{-1}u\}(x) \Big| dx}{ \int_{\Omega}|u(x)|dx} \\
  &=& \| V_- ( H_0 + \lambda)^{-1}\|_{L^1 \to L^1} 
\end{eqnarray*}
and
\[
\|  (H_0 + \lambda)^{-1}V_-\|_{L^{\infty}_{V_-} \to L^{\infty}_{V_-}} = \| (H_0+ \lambda)^{-1}V_-\|_{L^{\infty} \to L^{\infty}} =
\| V_- (H_0 + \lambda)^{-1}\|_{L^1 \to L^1} .
\]
By the Stein interpolation theorem we then obtain
\begin{eqnarray*}
\|  ( H_0 + \lambda)^{-1}V_-\|_{L^2_{V_-} \to L^2_{V_-}}&\leq& \|  ( H_0 + \lambda)^{-1}V_-\|_{L^1_{V_-} \to L^1_{V_-}}^{1/2} 
\| (H_0 + \lambda)^{-1}V_-\|_{L^{\infty}_{V_-} \to L^{\infty}_{V_-}}^{1/2} \\
&=&  \| V_- (H_0 + \lambda)^{-1}\|_{L^1 \to L^1}  ,
\end{eqnarray*}
so $V_-$ has zero operator bound with respect to $H_0$.
Applying \cite[Lemma 4.20]{dops} we conclude that Hypothesis (H) is satisfied. Let us note here that condition (\ref{l1}) is also related to certain integral conditions on $V$; see also \cite{dh,zy}.

\

\noindent
{\bf Example 2.} Suppose $V_-$ satisfies the weak Miyadera condition with respect to $H_0$: for any $\epsilon>0$ there exists $\delta >0$ such that
\beq
\int_0^{\delta} \|V_- e^{-tH_0}u\|_{L^1} dt \leq \epsilon \|u\|_{L^1}
\label{wm}
\eeq
for all $u\in L^1(\Omega)\cap L^2(\Omega)$. It is known \cite{mi,vo} that condition (\ref{l1}) is then satisfied, hence
Hypothesis (H) is satisfied.

\

Our first theorem reads:
\begin{theorem}
Let $2m>n$. Let $V$ be a real potential satisfying Hypothesis (H). Assume that the principal symbol $A(x,\xi)$ is strongly convex and that the principal coefficients
$a_{\ab}$, $|\alpha|=|\beta|=m$, belong in $W^{m-1,\infty}(\Omega)$. Then for any $\epsilon>0$ and $M>0$ there exists a constant $\Gamma_{\epsilon,M}$ such that the heat kernel of $H$ satisfies
\[
| K(t,x,y) |< \Gamma_{\epsilon,M}t^{-\frac{n}{2m}}\exp\Big\{ -(\sigma_m -\epsilon)
\frac{  d_M(x,y)^{\frac{2m}{2m-1}}}{   t^{\frac{1}{2m-1}}   } + \Gamma_{\epsilon, M}t   \Big\},
\]
for all $t>0$ and $x,y\in\Omega$.
\label{thm1}
\end{theorem}

Under additional assumptions we can obtain a better estimate that involves the actual Finsler distance $d(x,y)$ defined by
(\ref{finslermetric}) rather than the distances $d_M(x,y)$:
\begin{coro}
In addition to the assumptions of Theorem \ref{thm1} assume that $\ia$ $\Omega$ is bounded with $C^{m+1}$ boundary or $\Omega=\R^n$ and
$\ib$ the coefficients $a_{\ab}$ belong in $C^{m+1}(\Omega)$ and have bounded all derivatives of order up to $m+1$.
Then for any $\epsilon>0$ there exists $\Gamma_{\epsilon}$ such that
\be
| K(t,x,y) |< \Gamma_{\epsilon}t^{-\frac{n}{2m}}\exp\Big\{ -(\sigma_m -\epsilon)
\frac{  d(x,y)^{\frac{2m}{2m-1}}}{   t^{\frac{1}{2m-1}}   } + \Gamma_{\epsilon}t   \Big\},
\la{eq:114}
\ee
for all $t>0$ and $x,y\in\Omega$.
\label{coro}
\end{coro}
{\em Proof of Corollary \ref{coro}.} It is proved in \cite[Proposition 8 and Example p.595]{barbatis1999} that under the assumptions of the corollary
there holds
\begin{equation}
\frac{d_M(x,y)}{d(x,y)} \longrightarrow  1,  \qquad \mbox{ as }M\to +\infty ,
\label{fins_conv}
\end{equation}
uniformly in $x,y\in \Omega$. Estimate (\ref{eq:114}) then follows directly from
Theorem \ref{thm1} and  (\ref{fins_conv}). $\hfill\Box$

\

We next state a variation of Theorem \ref{thm1} which applies to
a wider class of operators. Let $D_A$ denote the distance in $L^{\infty}(\Omega)$ of the symbol
$A(x,\xi)$ to the class of all strongly convex symbols with coefficients in $W^{m-1,\infty}(\Omega)$; more precisely,
\[
D_A =\inf \max_{\subab} \| a_{\ab} - \hat{a}_{\ab}\|_{L^{\infty}(\Omega)},
\]
where the infimum is taken over all coefficient matrices $\{\hat{a}_{\ab}\}$ whose entries belong in $W^{m-1,\infty}(\Omega)$ and for which the symbol $\hat{A}(x,\xi)=\sum \hat{a}_{\ab}\xi^{\alpha+\beta}$ is strongly convex;
in particular $D_A=0$ if the symbol is strongly convex and the principal coefficients are uniformly continuous.
We then have
\begin{theorem}
Let $2m>n$. Let $V$ be a real potential satisfying Hypothesis (H).
For any $\epsilon>0$ and $M>0$ the heat kernel estimate
\[
| K(t,x,y) |< \Gamma_{\epsilon,M}t^{-\frac{n}{2m}}\exp\Big\{ -(\sigma_m -cD_A-\epsilon)
\frac{  d_M(x,y)^{\frac{2m}{2m-1}}}{   t^{\frac{1}{2m-1}}   } + \Gamma_{\epsilon, M}t   \Big\},
\]
is valid for some $\Gamma_{\epsilon,M}$ and all $t >0$ and $x,y\in\Omega$; here $c$ is a positive constant that depends on the operator
$H$ but not on $\epsilon$ or $M$.
\label{thm2}
\end{theorem}

\section{Proofs of Theorems}

Throughout this section we assume that $H$ is an operator defined via of the quadratic form (\ref{q}) where $V$ is a potential
satisfying Hypothesis (H). We do not yet assume that the the coefficients belong in $W^{m-1,\infty}(\Omega)$ or that the symbol
$A(x,\xi)$ is strongly convex; these assumptions will only be made when we arrive at equation (\ref{123}).

Our approach is based on Davies' exponential perturbation technique.
For any $M>0$ we define
\[
\cE_M= \{\phi\in C^m(\Omega) \, : \; |D^{\alpha}\phi(x)|\leq M  , \; x\in\Omega \, , \;\; 1 \leq |\alpha|\leq m\}.
\]
Let $\phi\in\cE_M$ be fixed. We define a sesquilinear form $Q_{\phi}$ by $\dom(Q_{\phi})=\dom(Q)$ and
\[
Q_{\phi}(u,v)= Q(e^{\phi}u,e^{-\phi}v) ;
\]
here $Q(\cdot,\cdot)$ denotes the sesquilinear form associated with the quadratic form $Q(\cdot)$.
We denote by $Q_{\phi}(\cdot)$ the quadratic form corresponding to the sesquilinear form $Q_{\phi}(\cdot,\cdot)$.
Let $H_{\phi}$ be the (non-self adjoint) operator associated to the form $Q_{\phi}(\cdot,\cdot)$, so that
$H_{\phi}=e^{-\phi}He^{\phi}$.
This conjugation induces canonically a functional calculus for $H_{\phi}$ via the relation $f(H_{\phi})=e^{-\phi}f(H)e^{\phi}$. In particular
$H_{\phi}$ is the generator of a strongly continuous semigroup given by
\begin{equation}
e^{-H_{\phi}t}=e^{-\phi}e^{-Ht} e^{\phi} .
\label{conj}
\end{equation}
\begin{lemma}
Assume that $2m>n$. Let $\phi\in\cE_M$ be fixed and let $k\in\R$  be such that
\[
{\rm Re}\, Q_{\phi}(u)  \geq -k \|u\|_{L^2}^2 \; , \qquad \mbox{ all }u\in  C^{\infty}_c(\Omega) .
\]
Then the heat kernel of $H$ satisfies
\beq
|K(t,x,y)|  \leq c_{\delta,M}t^{-\frac{N}{2m}}\exp\big\{\phi(y)-\phi(x) + (1+\delta)kt\big\}
\label{ult}
\eeq
for any $\delta>0$, all $t>0$ and $x,y\in\Omega$ and some constant $c_{\delta,M}$ which depends only on $\delta$ and $M$.
\label{lemma}
\end{lemma}
{\em Proof.} Let $Q_{0,\phi}(\cdot)$ denote the quadratic form defined as above for the free operator $H_0$ (rather than $H$).
The difference $Q_{0,\phi}(\cdot)-Q_0(\cdot)$ is of order smaller than $2m$ and this yields (see also \cite[Lemma 2]{davies1})
\[
|Q_{0,\phi}(u)-Q_0(u)|   \leq \epsilon Q_0(u) +c_{\epsilon,M} \|u\|^2_{L^2} ,
\]
for any $\epsilon>0$ and all $u\in C^{\infty}_c(\Omega)$. Hence we have from (\ref{q_0q})
\begin{eqnarray*}
|Q_{\phi}(u)-Q(u)|  &=& |Q_{0,\phi}(u)-Q_0(u)| \\
&\leq &\epsilon Q_0(u) +c_{\epsilon,M} \|u\|_{L^2}^2 \\
&\leq & \frac{\epsilon}{1-\epsilon} Q(u) + \Big( \frac{\epsilon c_{\epsilon}}{1-\epsilon}  +c_{\epsilon,M} \Big)
\|u\|_{L^2}^2 \, ,
\end{eqnarray*}
and therefore
\beq
Q(u) \leq  2 {\rm Re \,}Q_{\phi}(u) + c_M\|u\|_{L^2}^2   \; , \quad \mbox{ all } u\in C^{\infty}_c(\Omega)  .
\label{111}
\eeq
Now, let $u\in L^2(\Omega)$ be given and for $t>0$ let $u_t=e^{-H_{\phi}t}u$. By the multiplicative Sobolev inequality \cite[Lemma 16]{davies1} and by inequalities (\ref{q_0q}) and (\ref{111}) we have
\begin{eqnarray}
\|u_t\|_{L^{\infty}}&\leq& c Q_0(u_t)^{\frac{n}{4m}}\|u_t\|_{L^{2}}^{1-\frac{n}{2m}}\nonumber \\
&\leq& c \Big( Q(u_t) + \|u_t\|_{L^2}^2 \Big)^{\frac{n}{4m}}\|u_t\|_{L^{2}}^{1-\frac{n}{2m}}\nonumber \\
&\leq& c_M\Big({\rm Re}\, Q_{\phi}(u_t)  + \|u_t\|_{L^2}^2 \Big)^{\frac{n}{4m}} \|u_t\|_{L^2}^{1-\frac{n}{2m}} \nonumber\\
&\leq& c_M \Big( \| H_{\phi}u_t\|_{L^2}  \|u_t\|_{L^2} + \|u_t\|_{L^2}^2 \Big)^{\frac{n}{4m}}  \|u_t\|_{L^2}^{1-\frac{n}{2m}} \nonumber\\
&\leq& c_M \Big( \| H_{\phi}u_t\|_{L^2}^{\frac{n}{4m}}  \|u_t\|_{L^2}^{1-\frac{n}{4m}} + \|u_t\|_{L^2} \Big).  \label{2i}
\end{eqnarray}
Now, it follows from (\ref{q_0q}) and (\ref{111}) that for any $\delta>0$ there exists $c_{\delta}$ such that 
\[
\|u_t\|_{L^2} \leq e^{kt} \|u\|_{L^2}\;\; , \qquad \|H_{\phi}u_t\|_{L^2} \leq 
\frac{c_{\delta}}{t}e^{(k+\delta)t} \|u\|_{L^2} \, , \qquad t>0\, .
\]
This has been proved in \cite[Lemma 2.1]{bd} in the case $V=0$; since the proof in our case is identical, we omit further details.

Renaming $n\delta/4m$ as $\delta$ it follows from (\ref{2i}) that
\[
\|e^{-H_{\phi}t}u\|_{L^{\infty}}  \leq c_{\delta,M}e^{kt}\Big( t^{-\frac{n}{4m}}e^{\delta t}  +1  \Big) \|u\|_{L^2}
\leq c_{\delta,M}' t^{-\frac{n}{4m}}e^{(k+\delta) t}  \|u\|_{L^2}.
\]
Using duality we conclude that the semigroup $e^{-H_{\phi}t}$ maps $L^1\cap L^2$ into $L^{\infty}$ and
\[
\| e^{-H_{\phi}t} \|_{L^1\to L^{\infty}} \leq c_{\delta,M}t^{-\frac{n}{2m}}e^{(k+\delta)t}  .
\]
This together with (\ref{conj}) implies (\ref{ult}). $\hfill\Box$

\

\noindent
{\bf\em Proof of Theorem \ref{thm1}.} We shall now make use of the assumptions that $a_{\ab}\in W^{m-1,\infty}(\Omega)$
and that $A(x,\xi)$ is strongly convex. Let
\[
k_m =\Big(\sin(\frac{\pi}{4m-2})\Big)^{1-2m} \, .
\]
It has been proved in \cite[Proposition 6 and Lemma 7]{b2001} that for any $\epsilon , M >0$ there exists a constant $c_{\epsilon ,M}$ such that
\beq
{\rm Re}\, Q_{0,\lambda\phi}(u)\geq - \Big(\lambda^{2m}(k_m +\epsilon)  +c_{\epsilon,M}\Big)\|u\|_{L^2}^2,
\label{123}
\eeq
for all $\phi\in \cE_{A,M}$, all $\lambda>0$ and all $u\in C^{\infty}_c(\Omega)$ (the constant $c_{\epsilon,M}$ also depends on \nl
$\max_{\alpha,\beta} \max_{0\leq k\leq m-1}\|\nabla^k a_{\ab}\|_{L^{\infty}}$).
Moreover, using (\ref{ma}) and recalling (\ref{111}) (for $H_0$ rather than $H$) we obtain
\begin{eqnarray}
{\rm Re}\, Q_{\lambda\phi}(u) &=& {\rm Re}\, Q_{0,\lambda\phi}(u) +\int_{\Omega}V|u|^2dx \nonumber\\
&\geq&  {\rm Re}\, Q_{0,\lambda\phi}(u) -\epsilon Q_0(u) -c_{\epsilon}\|u\|_{L^2}^2 \nonumber \\
&\geq& (1-2\epsilon) {\rm Re}\, Q_{0,\lambda\phi}(u) -(\epsilon c_M +c_{\epsilon})\|u\|_{L^2}^2 .\label{stn}
\end{eqnarray}
From (\ref{123}) and (\ref{stn}) follows that for all $\epsilon>0$ small enough and for any $M>0$ there exists $c_{\epsilon,M}$ such that
\beq
{\rm Re}\, Q_{\lambda\phi}(u) \geq  -\Big(  \lambda^{2m}(k_m+\epsilon) +c_{\epsilon,M}\Big) \|u\|_{L^2}^2 \, ,
\qquad u\in C^{\infty}_c(\Omega).
\label{end}
\eeq
We complete the standard argument by first applying Lemma \ref{lemma} and then optimizing over all $\phi\in\cE_{A,M}$ and all
$\lambda>0$. Noting that
\[
\inf_{\lambda>0} \Big( -\lambda d_M(x,y) + \lambda^{2m}k_mt   \Big) =-\sigma_m\frac{d_M(x,y)^{\frac{2m}{2m-1}}}{t^{\frac{1}{2m-1}}}
\]
completes the proof of the theorem. $\hfill\Box$

\

\noindent
{\bf\em Proof of Theorem \ref{thm2}.} The main idea in the proof is that estimate (\ref{end}) is stable under perturbations of that are small in the $L^{\infty}$ norm: a perturbation of order $\delta$ (in $L^{\infty}$) results to a perturbation of order $\delta$ on the lower bound $k$ of Lemma~\ref{lemma}, and this results to a perturbation of order $\delta$ of the coefficient of $d_M^{2m/(2m1)}t^{-1/(2m-1)}$ in the Gaussian estimate. To see this, suppose $\hat{H}$ is an operator with coefficients $\{\hat{a}_{\ab}\}$ for which
\beq
{\rm Re}\, \hat{Q}_{\lambda\phi}(u) \geq  -\Big(  \lambda^{2m}(k_m+\epsilon) +c_{\epsilon,M}\Big) \|u\|_{L^2}^2 \, ,
\qquad u\in C^{\infty}_c(\Omega).
\label{end1}
\eeq
for all $\phi\in \cE_{\hat{A},M}$ and all $\epsilon,M>0$. Suppose now that $H$ is another operator satisfying the assumptions in Section \ref{section:hke} and such that $\|a_{\ab}-\hat{a}_{\ab}\|_{L^{\infty}}<\delta$, $|\alpha|=|\beta|=m$, where $\delta>0$ is small.
For any $u\in C^{\infty}_c(\Omega)$ we then have
\[
\Big|{\rm Re} Q_{\lf}(u) -{\rm Re} \hat{Q}_{\lf}(u)\Big| \leq c\delta \{ Q(u) + \lambda^{2m}\|u\|_{L^2}^2 \} +
c_M\delta (1+\lambda^{2m-1}) \|u\|_{L^2}^2
\]
This has been proved in detail in (cf. also \cite[eqn. (18)]{b2001}). The fact that the coefficient of $\lambda^{2m}$ can be estimated independently of $M$ is due to the fact that when $Q_{\lf}(u)$ and $\hat{Q}_{\lf}(u)$ are expanded into a polynomial of $\lambda$, the coefficient of $\lambda^{2m}$ involves only first-order derivatives of $\phi$ and not higher-order derivatives (see
also \cite[Lemma 3]{b2001}). Recalling also (\ref{111}) (with $\phi$ replaced by $\lambda\phi$) we thus obtain
\bean
{\rm Re \,} \hat{Q}_{\lf}(u)&\geq&  {\rm Re \,}Q_{\lf}(u) -c\delta \{ Q(u) + \lambda^{2m}\|u\|_{L^2}^2\} -
c_M\delta (1+\lambda^{2m-1}) \|u\|_{L^2}^2 \\
&\geq& (1-2c\delta) {\rm Re \,}Q_{\lf}(u)  -\bigg[  c\delta\lambda^{2m}+ c_M\delta (1+\lambda^{2m-1})\bigg] \|u\|_{L^2}^2 \\
&\geq& -\bigg[ (1-2c\delta)\Big(  \lambda^{2m}(k_m+\epsilon) +
c_{\epsilon,M}\Big) + c\delta\lambda^{2m}+ c_M\delta (1+\lambda^{2m-1})\bigg] \|u\|_{L^2}^2 .
\eean
We note that given $\epsilon_1>0$ the term in square brackets can be made smaller than
$\lambda^{2m}(k_m+c\delta +\epsilon_1) +c_{\epsilon_1,M}$, so estimate (\ref{end}) is true for the operator $H$ with $k_m$ being replaced by $k_m+c\delta$.

This leads to an estimate involving a constant $\sigma_m-c\delta$ and the distance $\hat{d}_M(x,y)$. To obtain an estimate with $d_M(x,y)$ we note that
there exists $c>0$ such that if $\hat{\phi}\in \cE_{\hat{A},M}$ then $(1+c\delta)^{-1}\phi \in \cE_{A,M}$. From this follows that
\bean
\hat{d}_M(x,y)&=&\sup \{\tilde{\phi}(y)- \tilde{\phi}(x)  : \tilde{\phi}\in \cE_{\hat{A},M} \} \\
&\geq& \sup \{ (1+c\delta)^{-1} \big( {\phi}(y)- {\phi}(x) \big)  : \phi\in \cE_{\hat{A},M} \} \\
&=& (1+c\delta)^{-1} d_{M}(x,y).
\eean
Combining the above concludes the proof of the theorem. $\hfill\Box$

\

\noindent
{\bf Acknowledgment} I thank J\"{u}rgen Voigt and Hendrik Vogt for useful comments. I also thank the referee whose suggestions
have improved the presentation of the article.

\end{document}